\documentclass[twoside, 12pt]{article}
\usepackage{amssymb, amsmath, mathrsfs, amsthm}
\usepackage{graphicx}
\usepackage{tikz}
\usepackage{color}
\usepackage[top=2cm, bottom=2cm, left=2cm, right=2cm]{geometry}
\usepackage{float, caption, subcaption}
\usepackage{ifpdf}
\usepackage{epstopdf}
\input{epsf}

\newcommand{\bd}{\begin{description}}
\newcommand{\ed}{\end{description}}
\newcommand{\bi}{\begin{itemize}}
\newcommand{\ei}{\end{itemize}}
\newcommand{\be}{\begin{enumerate}}
\newcommand{\ee}{\end{enumerate}}
\newcommand{\beq}{\begin{equation}}
\newcommand{\eeq}{\end{equation}}
\newcommand{\beqs}{\begin{eqnarray*}}
\newcommand{\eeqs}{\end{eqnarray*}}

\definecolor{DarkGreen}{rgb}{0.2, 0.6, 0.3}

\newtheorem{theorem}{Theorem}[section]

\newtheorem{lemma}{Lemma}[section]

\newtheorem{case}{Case}
\newtheorem{subcase}{Subcase}[case]
\newtheorem{claim}{Claim}

\newtheorem{proposition}{Proposition}[section]

\newtheorem{example}{Example}

\begin{document}
\title{\textbf{Asymptotic Bounds for CO-irredundant and Irredundant Ramsey Numbers} \footnote{Supported by the National
Science Foundation of China (Nos. 12061059, 11601254, and 11551001), Doctoral research project of Tianjin Normal University (52XB2111)
and the Qinghai Key Laboratory of Internet of Things Project
(2017-ZJ-Y21).} }

\author{Meng Ji\footnote{College of Mathematical Science, Tianjin Normal University, Tianjin, China {\tt
mji@tjnu.edu.cn}}, \ \ Yaping Mao\footnote{Faculty of Environment and Information Sciences, Yokohama
National University, 79-2 Tokiwadai, Hodogaya-ku, Yokohama 240-8501,
Japan. {\tt maoyaping@ymail.com}}, \ \ Ingo
Schiermeyer\footnote{Technische Universit{\"a}t Bergakademie
Freiberg, Institut f{\"u}r Diskrete Mathematik und Algebra, 09596
Freiberg, Germany. {\tt Ingo.Schiermeyer@tu-freiberg.de}}}
\date{}
\maketitle

\begin{abstract}
A set of vertices $X\subseteq V$ in a simple graph $G(V,E)$ is
irredundant (CO-irredundant) if each vertex $x\in X$ is either isolated in the
induced subgraph $G[X]$ or else has a private neighbor
$y\in V\setminus X$ ($y\in V$) that is adjacent to $x$ and to no other vertex
of $X$.
The irredundant Ramsey number $s(t_{1},\ldots,t_{l})$, CO-irredundant Ramsey number $s_{\operatorname{CO}}(t_{1},\ldots,t_{l})$, is the minimum $N$ such that every $l$-coloring of the edges of the complete graph $K_{N}$ on $N$ vertices has a monochromatic irredundant set, a monochromatic CO-irredundant set, of size $t_{i}$ for some $1\leq i\leq l$, respectively.
In this paper, firstly, we establish a lower bound for the irredundant Ramsey number $s(t_{1},\ldots,t_{l})$ by a random and probabilistic method. Secondly, we improve an upper bound for $s(3,9)$ such that $24\leq s(3,9)\leq 26$. Thirdly, using Krivelevich's lemma, we establish an asymptotic lower bound for the $\operatorname{CO}$-irredundant Ramsey number $s_{\operatorname{CO}}(m,n)$.\\[2mm]
\end{abstract}

\section{Introduction}

Let $G$ be a graph with vertex set $V(G)$ and edge set $E(G)$. We denote by $\delta(G)$ and
$\Delta(G)$ the minimum and maximum degrees of the vertices of $G$, respectively. For any subset $X\subseteq V(G)$, let $G[X]$ denote the subgraph induced
by $X$. Similarly, for any subset
$F \subseteq E(G)$, let $G[F]$ denote the subgraph induced by $F$.
A {\it path}\index{path} on $n$ vertices is denoted by
$P_n$, and a {\it cycle}\index{cycle} on $n$ vertices is denoted by
$C_n$. The {\it degree}\index{degree} of a vertex $v$ in a graph $G$,
denoted by $d_G(v)$, is the number of edges of $G$ incident with
$v$.
A graph $G$ is called {\it $k$-regular} if $d_G(v)=k$ for every $v\in V(G)$. The {\it join}\index{join} $G\vee H$ of two
disjoint graphs $G$ and $H$ is the graph with vertex set $V(G)\cup V(H)$
and edge set $E(G)\cup E(H)\cup \{uv\,|\,
u\in V(G), v\in V(H)\}$. The {\it union}\index{union} $G\cup H$ of two graphs $G$ and $H$ is the
graph with vertex set $V(G)\cup V(H)$ and edge set
$E(G)\cup E(H)$. For $V_{1},V_{2}\subset V$, we denote the number of edges between $V_{1}$ and $V_{2}$ by $e(V_{1},V_{2})$.
Any undefined concepts or notation can be found in \cite{BM}.

\subsection{(CO-irredundant) Irredundant Ramsey
numbers}
In 1978, Cockayne, Hedetniemi, and Miller \cite{Cockayne-Hedetniemi-Miller} introduced the concept of irredundance which is relevant for dominating sets.
A set of vertices $X\subseteq V$ in a simple graph $G(V,E)$ is
\emph{irredundant} if each vertex $x\in X$ is either isolated in the
induced subgraph $G[X]$ or else has a private neighbor
$y\in V\setminus X$ that is adjacent to $x$ and to no other vertex
of $X$.
 Farley and Schacham \cite{Farley-Schacham} defined the CO-irredundant set: given a graph $G=(V,E)$, a vertex subset $X\subseteq V$ is called
\emph{CO-irredundant} if every vertex $v\in X$ either contains
no neighbors in $X$ or else has a private neighbor
$y\in V$ that is adjacent to $x$ and to no other vertex of $X$. 
The \emph{irredundant Ramsey number}  $s(t_{1},\ldots,t_{l})$ (resp., \emph{CO-irredundant Ramsey number} $s_{\operatorname{CO}}(t_{1},\ldots,t_{l})$), is the minimum $N$ such that every $l$-coloring of the edges of the complete graph $K_{N}$ has a monochromatic $t_{i}$-element irredundant set (resp., $t_{i}$-element CO-irredundant set, say $\operatorname{COIR}_{t_{i}}$) for certain $1\leq i\leq l$. If $t_{1}=t_{2}=\ldots=t_{l}$, then we denote it by $s(t;l)$ (resp., $s_{\operatorname{CO}}(t;l))$). The definition of the \emph{Ramsey number} $r(t_{1},\ldots,t_{l})$ differs from $s(t_{1},\ldots,t_{l})$ in that the $t_{i}$-element irredundant set is replaced by a $t_{i}$-element independent set.
The \emph{mixed Ramsey number} $t(m,
n)$ is the smallest $N$ for which every red-blue coloring of the
edges of $K_N$ yields an $m$-element irredundant set in the blue
subgraph or an $n$-element independent set in the red subgraph.
Note that each independent set is an irredundant set, and each irredundant set is a $\operatorname{CO}$-irredundant set. Consequently, it follows that
$$
s_{\operatorname{CO}}(t_{1},\ldots,t_{l})\leq s(t_{1},\ldots,t_{l})\leq r(t_{1},\ldots,t_{l})
$$
and
$$
s_{\operatorname{CO}}(m,n)\leq s(m, n)\leq t(m,n)\leq r(m, n).
$$

The difficulty of obtaining exact values for irredundant Ramsey numbers is evidently comparable to that of obtaining exact values for classical Ramsey numbers.
Brewster, Cockayne, and Mynhardt \cite{BrewsterCockayneMynhardt}
proposed irredundant Ramsey numbers and established the values $s(3,3) =
6$, $s(3,4) = 8$, $s(3,5) = 12$ and $s(3,6) = 15$ was established
in \cite{BrewsterCockayneMynhardtII}. It
was furthermore shown in \cite{Hattingh} that $18\leq s(3,7)\leq
19$ and Chen and Rousseau proved that
$s(3,7) = 18$ in \cite{ChenRousseau}, and Cockayne et al. in \cite{CockayneExooHattinghMynhardt} obtained that $s(4, 4) = 13$.
The values $t(3,3)=6, t(3,4)=9, t(3,5)=12$ and $t(3,6)=15$ have been shown in \cite{CockayneHattinghKokMynhardt, HenningOellermann}.
Burger, Hattingh, and Vuuren
\cite{BurgerHattinghvanVuuren} proved that $t(3,7)=18$ and
$t(3,8)=22$. Burger and Vuuren \cite{BurgerVuuren} obtained that
$s(3,8)=21$.
The following table lists all known Ramsey numbers $s(3,n), t(3,n)$ and $r(3,n)$ for $3 \leq n \leq 9.$ 

\begin{center}
\begin{center}
Table. 1: Exact known Ramsey numbers
\end{center}\label{fig4-3}
\begin{tabular}{c|c|c}
s(m,n) & t(m,n) & r(m,n)\\
\hline
s(3,3)=6 & t(3,3)=6 & r(3,3)=6\\
s(3,4)=8 & t(3,4)=9 & r(3,4)=9\\
s(3,5)=12 & t(3,5)=12 & r(3,5)=14\\
s(3,6)=15 & t(3,6)=15 & r(3,6)=18\\
s(3,7)=18 & t(3,7)=18 & r(3,7)=23\\
s(3,8)=21 & t(3,8)=22 & r(3,8)=28\\
- & - & r(3,9)=36\\
\end{tabular}
\end{center}

Chen, Hattingh, and Rousseau \cite{chen1993}, Erd\H{o}s and Hattingh \cite{ErdosHattingh}, and Krivelevich \cite{Krivelevich} have obtained several asymptotic bounds for irredundant  Ramsey numbers $s(m,n)$ and mixed Ramsey number $t(m,n)$. What's more, problems related to irredundant Tur\'an numbers has been studied in \cite{
CockayneHattinghKokMynhardt}. Furthermore, for $s_{\operatorname{CO}}(m,n)$, several exact values were given by Cockayne, MacGillivray and Simmons in \cite{Cockayne-MacGillivray-Simmons}. However, the asymptotic bounds for $s_{\operatorname{CO}}(m,n)$ are not given.
In 1994, Cockayne and Mynhardt \cite{CockayneMynhardt} gave the exact value $s(3,3,3)=13$.

For the $2$-coloring of the edges of $K_{N}$, we call it red-blue coloring and we call two kinds of monochromatic edge-induced subgraphs the red graph $\langle R\rangle$ and the blue graph $\langle B\rangle$.
If $Y$ is an $m$-element irredundant set in $\langle B\rangle$, then for some $k\leq m$, there exist $k$ vertices of $Y$ that have private neighbors in $\langle B\rangle$ and the remaining $m-k$ vertices of $Y$ in the induced subgraph of the $\langle B\rangle$ are isolated. With the $k$ vertices in $Y$ and their private neighbors, there is a $m+k$-element set in which all but $2{k\choose 2}$ of the ${m+k\choose 2}$ internal edges are completely
determined. So $\langle R\rangle$ contains one or more of the graphs from the graph family $\{K_{m},K_{m-k}+(K_{k,k}-kK_{2}),K_{m,m}-mK_{2}|3\leq k\leq m-1\}$ where the graph $K_{k,k}-kK_{2}$ is obtained by removing $k$ independent edges from the complete bipartite graph $K_{k,k}$. Clearly, $K_{m-k}+(K_{k,k}-kK_{2})$ has $m+k$ vertices and ${m+k\choose 2}-2{k\choose 2}-k$ edges.

\subsection{Our results}
Sawin \cite{WillSawin} proved a lower bound for $r(t;l)$ using the random mathod. We will prove a similar result for $s(t;l)$ in Section \ref{Multicolor}.

\begin{theorem}\label{multi-irredundant}
For $l> 2$ and sufficiently large $t$, we have
$$
s(t;l)\geq \left(\frac{t}{3}\right)^{\frac{3-l}{2}}2^{\frac{lt-t}{2}}.
$$
\end{theorem}

Burger and Vuuren \cite{BurgerVuuren} showed the upper and lower
bounds for $s(3,9)$ and $t(3,9)$.

\begin{theorem}[{\upshape{\bf Burger and Vuuren} \cite{BurgerVuuren}}]\label{lem4}
$24\leq s(3,9)\leq t(3,9)\leq 27$.
\end{theorem}

In Section \ref{upper-bound-$s(3,9)$}, we derive an upper bound for $s(3,9)$.
\begin{theorem}\label{thm2}
$24\leq s(3,9)\leq26$.
\end{theorem}
In Section \ref{section1}, we establish a lower bound for $s_{\operatorname{CO}}(m,n)$.
\begin{theorem}\label{CO-irredund-thm}
For each $m\geq 3$, there is a positive constant $c_{m}$ such that
\begin{equation*}
s_{\operatorname{CO}}(m,n)>c_{m}\left(\frac{n}{\log n}\right)^{\rho(\mathcal{H})},
\end{equation*}
where $\mathcal{H}$ is a graph family of $((K_{|W_{1}|}-tK_{2})\vee(K_{|Y_{1}|,|Y_{3}|}-|Y_{1}|K_{2}))\vee(K_{|W_{2}|,|Y_{2}|}\{Y_{2}\}-|W_{2}|K_{2})$ for $t\leq |W_{1}|\leq m-2$.
\end{theorem}

\section{A lower bound for $s(t;l)$}\label{Multicolor}

In this section, we obtain a lower bound for the irredundant Ramsey number $s(t;l)$. 

\noindent {\bf Proof of Theorem \ref{multi-irredundant}}
Let $b_{m}$ be the number of complete subgraphs of order $m$ in a graph $G$. Then a function was defined by the formulation
$$
f(n,s,t)=\min\{b_{s}(G):|G|=n,\,b_{t}(\bar{G})=0\},
$$
where $s\geq 3$ and $t\geq 3$.
Namely, it represents the minimum number of independent sets of size $s$ in a graph $G$ with $n$ vertices which contains no clique of size $t$, and is to study the number of complete subgraphs contained in a given graph.


Let $c_{s,t}=\lim_{n\rightarrow\infty}f(n,s,t)/{n\choose s}$, that is, to be the infimum, over graphs $G$ with
no $t$-clique, of the probability that $\{v_{1}, \ldots, v_{s}\}$ is an
independent set for the vertices $v_{1}, \ldots, v_{s}$ of
$G$ chosen independently and uniformly at random.
Nikiforov \cite{Nikiforov} first applied it to classical Ramsey number problem.

Clearly, we could use the method to obtain a new lower bound for irredundant Ramsey number $s(t;l)$. Let $g(n,s,t)$ be the minimum number of $s$-element irredundant
sets in a graph $G$ with $n$ vertices that contains no $\bar{I_{t}}$ of order $t$, where $\bar{I_{t}}$ denotes the complement graph of a $t$-element irredundant set in $K_{n}$.
Let
$$
w_{s,t}=\lim_{n\rightarrow\infty}\frac{g(n,s,t)}{{n\choose s}\left(1+\sum_{k=3}^{s}{s\choose k}{n-s\choose k}\right)}
$$
be the infimum over graphs $G$ with no $\bar{I_{t}}$ of order $t$ of
the probability that $\{v_{1}, \ldots, v_{s}\}$ is an irredundant set for the vertices $v_{1}, \ldots, v_{s}$ of $G$ chosen independently and
uniformly at random.

By proving the following Lemma \ref{lem1} and Lemma \ref{lem-1}, we immediately have Theorem \ref{multi-irredundant}. Before that, we first show two propositions for convenience of the proofs of those Lemmas.

\begin{proposition}\label{proposition1}
Let $n,t,k$ be three positive integers such that $n=\lceil 2^{t/2}\sqrt{t/3}\rceil$ and $k\geq 3$. Let
$$
b_{k}={t\choose k}{n-t\choose k}k!2^{{t+k\choose 2}-2{k\choose 2}}
$$
for $3\leq k\leq t$. For all sufficiently large value $t$, we have
$\max\{b_{k}\,|\,3\leq k\leq t\}=b_3$ or $b_t$.
\end{proposition}
\begin{proof}
For $3\leq k\leq(1-\epsilon)t/2$, since $n=\lceil 2^{t/2}\sqrt{t/3}\rceil$, it follows that
\begin{equation*}
\begin{split}
\frac{b_{k+1}}{b_{k}}&=\frac{(t-k)(n-t-k)}{k+1}2^{k-t}\\[0.2cm]
&\leq (t-k)(n-t-k)2^{k-t-2}~(since~k\geq 3)\\[0.2cm]
&\leq tn2^{k-t-2}\\[0.2cm]
&\leq t 2^{k-t-2}(2^{t/2}\sqrt{t/3}+1)\\[0.2cm]
&\leq t 2^{k-1-t/2}\sqrt{t/3}<1.
\end{split}
\end{equation*}
This means that the sequence $b_{3},\ldots,b_{t}$ is decreasing for $3\leq k\leq(1-\epsilon)t/2$.

If $t/2< k< t$, then
\begin{equation*}
\begin{split}
\frac{b_{k+1}}{b_{k}}&=\frac{(t-k)(n-t-k)}{k+1}2^{k-t}\\[0.2cm]
&\geq \frac{(t-k)(2^{t/2}\sqrt{t/3}-t-k)}{k+1}2^{k-t}\\[0.2cm]
&\geq \frac{(t-k)(2^{t/2}\sqrt{t/4})}{k+1}2^{k-t}\\[0.2cm]
&\geq \frac{\sqrt{t}(t-k)}{k+1}2^{k-1-t/2},
\end{split}
\end{equation*}
where the second inequality holds, since~the~sufficiently~large~$t$.

If $k=\lceil t/2\rceil$ or $k=t-1$, then
$$\frac{b_{\lceil t/2\rceil+1}}{b_{\lceil t/2\rceil}}> \frac{\lfloor t/2\rfloor\sqrt{t}}{2\left(\lceil t/2\rceil+1\right)}>1
$$
or
$$
\frac{b_{t}}{b_{t-1}}\geq\frac{2^{t/2-2}}{\sqrt{t}}> 1.
$$
If $t/2<k\leq t-3$ and $t$ is sufficiently~large, then
\begin{equation*}
\begin{split}
\frac{b_{k-1}b_{k+1}}{b_{k}^{2}}&=2\left(\frac{t-k}{t-k+1}\right)\left(\frac{n-t-k}
{n-t-k+1}\right)\left(\frac{k}{k+1}\right)\\[0.2cm]
&> 2\left(1-\frac{1}{1+t/2}\right)\left(1-\frac{1}{2^{t/2}\sqrt{t/4}}\right)\left(1-\frac{2}{t}\right)\\[0.2cm]
&>1.
\end{split}
\end{equation*}
This means that the sequence $b_{3},\ldots,b_{t}$ is increasing for $t/2< k< t$.

So, for all sufficiently large value $t$, we have
$\max\{b_{k}\,|\,3\leq k\leq t\}=b_3$ or $b_t$.
\end{proof}

\begin{proposition}\label{proposition2}
Let $n,t,k,h$ be positive integers such that $n=\lceil 2^{t/2}\sqrt{t/3}\rceil$ and $k\geq 3$. Let
$$
a_{k}={h\choose k}{n-h\choose k}k!2^{{h+k\choose 2}-2{k\choose 2}}
$$
for $3\leq k\leq h\leq t$. For all sufficiently large value $t$, we have
$\max\{a_{k}\,|\,3\leq k\leq t\}=a_3$ or $a_t$.
\end{proposition}

\begin{proof}
We have
$$
\frac{a_{k+1}}{a_{k}}=\frac{(h-k)(n-h-k)}{k+1}2^{k-h}.
$$
If $h>t/2$, then
$$
\frac{a_{k+1}}{a_{k}}<\left(\frac{hn}{4}\right)2^{k-h}=\frac{h}{4}\left(\frac{t}{3}\right)^{1/2}2^{k-h+t/2}\leq \frac{h}{4}\left(\frac{t}{3}\right)^{1/2}2^{\epsilon t}<1
$$
for $3\leq k\leq (1-\epsilon)(h-t/2)$ and for the sufficiently large $t$. And for $h-t/2<k<h$, it follows that
\begin{equation*}
\begin{split}
\frac{a_{k+1}}{a_{k}}&=\frac{(h-k)(n-h-k)}{k+1}2^{k-h}\\[0.2cm]
&\geq \frac{(h-k)(2^{t/2}\sqrt{t/3}-t-k)}{k+1}2^{k-h}\\[0.2cm]
&\geq \frac{(h-k)2^{t/2}\sqrt{t}/2}{k+1}2^{k-h}\\[0.2cm]
&\geq \frac{(h-k)\sqrt{t}/2}{k+1}2^{k-h+t/2}
\end{split}
\end{equation*}
for the sufficiently large $t$.

If $k=\lceil h-t/2\rceil+\mu$ for the positive integer $1\leq \mu\leq4$ or $k=h-1$, then it follows from $h>t/2$ that
\begin{equation*}
\begin{split}
\frac{a_{\lceil h-t/2\rceil+\mu+1}}{a_{\lceil h-t/2\rceil+\mu}}
&\geq \frac{(h-k)\sqrt{t}/2}{k+1}2^{k-h+t/2}\\[0.2cm]
&\geq \frac{(t/2-\mu-1)\sqrt{t}/2}{t+1}2^{\mu+1}\\[0.2cm]
&=O(t^{1/2})>1\\[0.2cm]
\end{split}
\end{equation*}
or
$$
\frac{a_{h}}{a_{h-1}} >\frac{\sqrt{t}2^{t/2-1}}{2h}> 1.
$$
If $h-t/2+5<k\leq h-2$, then
\begin{equation*}
\begin{split}
\frac{a_{k-1}a_{k+1}}{a_{k}^{2}}
&=2\left(1-\frac{1}{h-k+1}\right)\left(1-\frac{1}
{n-h-k+1}\right)\left(1-\frac{1}{k+1}\right)\\[0.2cm]
&>2\left(1-\frac{1}{3}\right)\left(1-\frac{1}
{1+2^{t/2}\sqrt{t}/2}\right)\left(1-\frac{1}{5+1}\right)\\[0.2cm]
&>\frac{10}{9}\left(1-\frac{1}
{1+2^{t/2}\sqrt{t}/2}\right)>1,\\[0.2cm]
\end{split}
\end{equation*}
for the sufficiently large $t$.

So, for $h>t/2$, the sequence $a_{3},\ldots,a_{h}$ is decreasing for $3\leq k\leq(1-\epsilon)(h-t/2)$ and is increasing for $h-t/2< k< h$ and the largest one must be $a_{3}$ or $a_{h}$.
\end{proof}

\begin{lemma}\label{lem1}
For sufficiently large $t$, we have
$$
w_{t,t}\leq \left(\frac{t}{3}\right)^{\frac{t}{2}}2^{\frac{-t^{2}}{2}+o(t^{2})}.
$$
\end{lemma}
\begin{proof}
Let $G$ be a random graph with $n$ vertices ($n$ will be chosen later), where each pair of vertices is connected by an edge with $p$. Let $v_{1},\ldots,v_{t}$ be uniformly distributed random variable in $[n]$, independent from each other and from $G$. For $w_{t,t}$ we have
$$
w_{t,t}\leq \frac{\mathbb{P}(\{v_{1},\ldots,v_{t}\}\, \mbox{is an irredundant set\,})}{\mathbb{P}(G \,\mbox{contains no}\, \bar{I}_{t})}.
$$
For the denominator, we first consider that $G \,\mbox{contains an}\, \bar{I}_{t}$.
\begin{equation}\label{equ1}
\mathbb{P}(G \,\mbox{contains an}\, \bar{I}_{t})\leq
{n\choose t}\left[p^{t\choose 2}+\sum_{k=3}^{t}{t\choose k}{{n-t}\choose k}k!p^{{{k+t}\choose 2}-2{k\choose 2}-k}(1-p)^{k}\right].
\end{equation}
We take $n=\sqrt{t/3}p^{-t/2}$ and $p=1/2$. From (\ref{equ1}) we have that
\begin{equation}\label{equ2}
\mathbb{P}(G \,\mbox{contains an}\, \bar{I}_{t})\leq
{n\choose t}\left[2^{-{t\choose 2}}+\sum_{k=3}^{t}{t\choose k}{{n-t}\choose k}k!2^{-{{{k+t}\choose 2}+2{k\choose 2}}}\right].
\end{equation}
Since
$b_{t}={{n-t}\choose t}t!2^{-t^{2}}\sim \left(\sqrt{t/3}\right)^{t}2^{-t^{2}/2}$
and
$b_{3}=\left(t\choose 3\right)\left({n-t}\choose t\right)3!2^{-{{t+3}\choose 2}+6}\sim 4/3\left(t\sqrt{t/3}\right)^{3}2^{-t}2^{-t^{2}/2}$
for the sufficiently large $t$. Then we have $b_{3}=o(b_{t})$ and $2^{-{t\choose 2}}=o(b_{t})$.
From (\ref{equ2}), we have
\begin{equation*}
\begin{split}
\mathbb{P}(G \,\mbox{contains an}\, \bar{I}_{t})&\leq
{n\choose t}\left[2^{-{t\choose 2}}+\sum_{k=3}^{t}{t\choose k}{{n-t}\choose k}k!2^{-{{{k+t}\choose 2}+2{k\choose 2}}}\right]\\[0.2cm]
&\leq {n\choose t}2t{{n-t}\choose t}t!2^{-t^{2}}(By~Proposition~\ref{proposition1})\\[0.2cm]
&\leq {n\choose t}2tn^{t}2^{-t^{2}}
\leq \left(\frac{ne}{t}\right)^{t}2tn^{t}2^{-t^{2}}\\[0.2cm]
&\leq 2t\left(\frac{e}{3}\right)^{t}=o(1),
\end{split}
\end{equation*}
which means that $\mathbb{P}(G \,\mbox{contains no}\, \bar{I}_{t})=1-o(1)$.

For the numerator taking $h$ to be the of size of $\{v_{1},\ldots,v_{t}\}$ and $p=1/2$, we have
\begin{equation*}
\begin{split}
&\mathbb{P}(\{v_{1},\ldots,v_{t}\}\, \mbox{is an irredundant set\,})\\[0.2cm]
&\leq  \sum_{h=1}^{t}\frac{{t \brace h}{n\choose h}h!\left((1/2)^{{h\choose 2}}+\sum_{k=3}^{h}{h\choose k}{{n-h}\choose k}k!(1/2)^{{{k+h}\choose 2}-2{k\choose 2}}\right)}{n^{t}+\sum_{k=3}^{t}n^{t+k}}\\[0.2cm]
&\leq \sum_{h=1}^{t}\frac{{t \brace h}\left((1/2)^{{h\choose 2}}+\sum_{k=3}^{h}{h\choose k}{{n-h}\choose k}k!(1/2)^{{{k+h}\choose 2}-2{k\choose 2}}\right)}{n^{t-h}+\sum_{k=3}^{t}n^{t+k-h}}\\[0.2cm]
&\leq \sum_{h=1}^{t}\frac{{t \brace h}\left((1/2)^{{t\choose 2}}+\sum_{k=3}^{t}{t\choose k}{{n-t}\choose k}k!(1/2)^{{{k+t}\choose 2}-2{k\choose 2}}\right)}{{tn^{t-h}}}\\[0.2cm]
&\leq t!\max_{0\leq h\leq t}\left\{\frac{tn^{t}2^{-t^{2}}}{{tn^{t-h}}}\right\}~~(By~\sum_{h=1}^{t}{t \brace h}\leq t!~and~Proposition~\ref{proposition2})\\[0.2cm]
&\leq 2^{t\log t}\max_{0\leq h\leq t}\left\{n^{h}2^{-t^{2}}\right\}~~\left(Since~t!\leq2^{t\log t}\right)\\[0.2cm]
&\leq 2^{t\log t}\max_{0\leq h\leq t}\left\{\left(\frac{t}{3}\right)^{\frac{t}{2}}2^{\frac{ht}{2}}2^{-t^{2}}\right\}
\leq 2^{t\log t}\left(\frac{t}{3}\right)^{\frac{t}{2}}2^{\frac{t^{2}}{2}}2^{-t^{2}}\\[0.2cm]
&\leq \left(\frac{t}{3}\right)^{\frac{t}{2}}2^{\frac{-t^{2}}{2}+o(t^{2})}.
\end{split}
\end{equation*}


\noindent{\bf Remark:} ${t \brace h}$ are the Stirling numbers of the second kind.

Combine the upper bound for $\mathbb{P}(\{v_{1},\ldots,v_{t}\}\, \mbox{is an irredundant set\,})$ with $\mathbb{P}(G \,\mbox{contains no}\, \bar{I}_{t})$, we conclude that
\begin{align*}
w_{t,t}&\leq \frac{\mathbb{P}(\{v_{1},\ldots,v_{t}\}\, \mbox{is an irredundant set\,})}{\mathbb{P}(G \,\mbox{contains no}\, \bar{I}_{t})}\\[0.2cm]
&\leq \frac{\left(\frac{t}{3}\right)^{\frac{t}{2}}2^{\frac{-t^{2}}{2}+o(t^{2})}}{1-o(1)}
= \left(\frac{t}{3}\right)^{\frac{t}{2}}2^{\frac{-t^{2}}{2}+o(t^{2})}.
\end{align*}
\end{proof}

\begin{lemma}\label{lem-1}
For $l> 3$ and sufficiently large $t$, we have
$$
s(t;l)\geq w_{t,t}^{-\frac{l-2}{t}}2^{\frac{t}{2}}\sqrt{\frac{t}{3}}.
$$
\end{lemma}
\begin{proof}
We first fix a graph $G$ satisfying two conditions:
\begin{itemize}
\item there is no $\bar{I}_{t}$ in $G$;

\item the probability that $t$ random vertices of $G$ form an irredundant set is at most
$w_{t,t}+\epsilon$ for $\epsilon \ (=\epsilon(l,t)) >0$.
\end{itemize}
Let $N=\left\lfloor w_{t,t}^{-\frac{l-2}{t}}2^{\frac{t}{2}}\sqrt{\frac{t}{3}}\right\rfloor$.
We construct an $l$-edge-coloring of $K_{N}$ without a monochromatic $\bar{I}_{t}$ as follows.
Let $f_{1},\ldots,f_{q-2}$ be $q-2$ uniform, independent and random functions 
from the vertex set of $K_{N}$ to the vertex set of $G$, all independent of one another, satisfying that,
for each subgraph $W\subset K_{N}$ with $|W|=2t$ and for $1\leq i\leq q-2$, $f_{i}$ maps each vertex of $V(W)$ into distinct vertices of a subset of $|W|$ vertices of $G$.

As above mapping rules, for $v_{1},v_{2}\in V(K_{N})$, if there are
some functions $f_{i}$ for $i\in[q-2]$ such that $f_{i}(v_{1}),f_{i}(v_{2})$ form an edge in $G$, 
then we define the color 
$\chi(v_{1}v_{2})=\min\{i\,|\,f_{i}(v_{1})f_{i}(v_{2})\in E(G)$ for $i\in[q-2]\}$, which means that we color the edge $v_{1}v_{2}$ by the minimum subscript of all the functions satisfying $f_{i}(v_{1})f_{i}(v_{2})\in E(G)$ for $i\in[q-2]$. 
If there is no some $i$ such that $f_{i}(v_{1}),f_{i}(v_{2})$ form an edge in $G$,
then we color randomly this edge $v_{1}v_{2}$ with $q-1$ or $q$, each with probability $1/2$, independently for each remaining edges.

Since $G$ contains no $\bar{I}_{t}$, it follows that there exists no set of $t$ vertices sent to an $\bar{I}_{t}$ by any $f_{i}$, $1\leq i\leq l-2$.
Thus, there exists no a monochromatic $\bar{I}_{t}$ for color $i$. So it suffices to prove that the
probability that there is an $\bar{I}_{t}$ for the remaining two colors is less than $1$. By the rule
of $l$-coloring, for $T$ to be an $\bar{I}_{t}$, it contains no edges of color $i$ for each $i\in [l-2]$, and hence the probability of occurrence is at most $(w_{t,t}+\epsilon)^{l-2}$. If, for example, $T$ is an
irredundant set with color $l-1$, then for some subset
$\{x_{1},\ldots,x_{k}\}\subset V(K_{N})\setminus T$, the private neighbors of
the vertices $y_{1},\ldots,y_{k}$ are not isolated in $T$. So there are exactly $2{k\choose 2}$
edges such that their colors are non-deterministic, and the coloring of the remaining ${t+k\choose 2}-2{k\choose 2}$
edges is completely determined.
  Let $X_{T}$ denote the indicator random variable that takes the value $i$ if $T$
is an $\bar{I}_{t}$ with a monochromatic color $l-1$ or color $l$. We denote the size of $T$ by $t$, and let $X$ be the random variable that counts the number of $\bar{I}_{t}$ of color $l-1$ or $l$.
For a fixed set $T$ with $k$ non-isolated vertices, the choices for the set of non-isolated vertices,
the private neighbors, and the matching between these two sets are ${t\choose k}$ and ${N-t\choose k}$, and $k!$, respectively. So we have that
$$
\mathbb{E}(X)=\sum_{T}E(X_{T})\leq 2{N\choose t}\left[2^{-{t\choose 2}}+\sum_{k=3}^{t}{t\choose k}{N-t\choose k}k!2^{{t+k\choose 2}-2{k\choose 2}}\right](w_{t,t}+\epsilon)^{l-2}.
$$
Let
$$
q_{k}={t\choose k}{N-t\choose k}k!2^{{t+k\choose 2}-2{k\choose 2}}
$$
for $3\leq k< t$. Naturally,

\begin{equation*}
\begin{split}
\frac{q_{k+1}}{q_{k}}&=\frac{(t-k)(N-t-k)}{k+1}2^{k-t}\\[0.2cm]
&\geq \frac{w_{t,t}^{-\frac{l-2}{t}}2^{\frac{t}{2}}\sqrt{\frac{t}{3}}-t-k-1}{t}2^{k-t}\\[0.2cm]
&\geq  \frac{\left(\left(\frac{t}{3}\right)^{\frac{t}{2}}2^{\frac{-t^{2}}{2}}\right)
^{-\frac{l-2}{t}}2^{\frac{t}{2}}\sqrt{\frac{t}{3}}-t-k-1}{t}2^{k-t}~~(By~Lemma~\ref{lem1})\\[0.2cm]
&\geq \frac{2^{k+t/2}\sqrt{3}-t-k-1}{t\sqrt{t}}~~(Since~l\geq 4)\\[0.2cm]
&> 1.
\end{split}
\end{equation*}

It is easily seen that the largest term of the sum is $q_{t}$.
Then $q_{t}={{n-t}\choose t}t!2^{-t^{2}}\sim \left(\sqrt{t/3}\right)^{t}2^{-t^{2}/2}$, $q_{3}=\left(t\choose 3\right)\left({n-t}\choose t\right)3!2^{-{{t+3}\choose 2}+6}\sim 4/3\left(t\sqrt{t/3}\right)^{3}2^{-t}2^{-t^{2}/2}$ for sufficiently large $t$ and $2^{-{t\choose 2}}=o(q_{t})$.
\begin{align*}
\mathbb{E}(X)=\sum_{T}E(X_{T})&\leq 2{N\choose t}\left[2^{-{t\choose 2}}+\sum_{k=3}^{t}{t\choose k}{N-t\choose k}k!2^{{t+k\choose 2}-2{k\choose 2}}\right](w_{t,t}+\epsilon)^{l-2}.\\[0.2cm]
&\leq 2t{N\choose t}{{N-t}\choose t}t!2^{-t^{2}}(w_{t,t}+\epsilon)^{l-2}\\[0.2cm]
&< 2t{N\choose t}N^{t}2^{-t^{2}}w_{t,t}^{l-2}<\left(\frac{Ne}{t}\right)^{t}2tN^{t}2^{-t^{2}}w_{t,t}^{l-2}\\[0.2cm]
&\leq 2t\left(\frac{e}{3}\right)^{t}<1.
\end{align*}
Clearly, the probability that there is no monochromatic $\bar{I}_{t}$ in $K_{N}$ is more than $0$.
\end{proof}

\section{An upper bound for $s(3,9)$}\label{upper-bound-$s(3,9)$}

An $s(m,n,p)$ (resp., $(r(m,n,p))$) coloring is a red-blue
edge-coloring $(R,B)$ of the complete graph $K_p$ which satisfies
neither of the two conditions of the irredundant Ramsey number
$s(m,n)$ (resp., $r(m,n)$), and hence showing that $s(m,n)>p$ (resp.,
Ramsey number $r(m,n)>p$).

\begin{lemma}{\upshape \cite{BurgerVuuren}}\label{zuidazuixiaodu}
For two integers $m,n\geq 2$ let $x\in\{s,t\}$ be a Ramsey
number. Then we have
\begin{itemize}
  \item $\Delta(R)<x(m-1,n)$ and $\Delta(B)<x(m,n-1)$ in any $x(m,n,p)$ coloring $(R,B)$.
  \item $\delta(R)\geq p-x(m,n-1)$ and $\delta(B)\geq p-x(m-1,n)$ in any $x(m,n,p)$ coloring $(R,B)$.
\end{itemize}
\end{lemma}

We will use the following two theorems in our proof.
\begin{theorem}[{\upshape {\bf Brewster, Cockayne, and Mynhardt} \cite{BrewsterCockayneMynhardt}}]\label{lem2}
The blue subgraph of a red-blue edge coloring of a complete graph
contains an irredundant set of cardinality $3$ if and only if the
red subgraph contains a $3$-cycle or an induced $6$-cycle.
\end{theorem}

\begin{theorem}[{\upshape{\bf Hattingh} \cite{Hattingh}}]\label{Hattingh}
Let $(R,B)$ be a red-blue coloring of the edges of a complete
graph $K_{N}$ in which $\langle B\rangle$ contains no $3$-element
irredundant set. For each vertex $v\in V(K_{N})$, we partition
$N_{B}(v)$ into two parts $D_{2}$ and $D_{>2}$, where each vertex $u$ of $N_{B}(v)$ belongs
to $D_{2}$ if the distance between $u$ and $v$ in $\langle R\rangle$ is two,
otherwise it belongs to $D_{>2}$. Let $X$ be any subset of $N_{B}(v)$
that contains at most one vertex from $D_{>2}$. Then $\langle
X\rangle_{R}$ is bipartite.
\end{theorem}

For the sake of simplicity, we define some notations. A graph $G$ is called an \emph{$(m,n)$-graph}, if it neither contains an
$m$-irredundant set in $G$ nor an $n$-irredundant set in
$\overline{G}$.  Denote the distance
in $G$ from $u$ to $v$ by $d(u,v)$; let
\begin{center}
$D_{i}=\{u\,|\,d(u,v)=i\}$ and $D_{>i}=\{u\,|\,d(u,v)>i\}$.
\end{center}
Let $(X,Y)$ be a bipartition of $G[D_{2}]$ and let $(X_{i},Y_{i})$
be bipartitions of the components of $G[D_{2}]$ for $i=1,2,\cdots,k$.
We denote the number of connected components of $(X,Y)$ by $C$. Without
loss of generality, assume that $D_{2}=X\cup Y$, $X=\cup X_{i}$,
$Y=\cup Y_{i}$ and $|X_{i}|\geq |Y_{i}|$ for $i=1,2,\cdots,C$. We say
there is an $UU'$-edge if $N(U)\cap U'\neq \emptyset$ for two disjoint vertex sets
$U,U'$ of $G$.\\[0.2mm]

\noindent {\bf Proof of Theorem \ref{thm2}.}
Suppose, to the contrary, that $G$ is a $(3,9)$-graph
of order $26$. Let $G$ be the red subgraph. Then it follows
\begin{center}
$5\leq\delta(G)\leq \Delta(G)<9$
\end{center}
by Lemma \ref{zuidazuixiaodu}. Let $v$ be a vertex of $G$ with
$d(v)=\Delta(G)$. Hence, $d(v)=5,6,7$, or $8$.
Since the order of $G$ is $26$ and $d(v)=\Delta(G)\leq 8$,
then $|D_{\geq2}|\geq 26-8-1=17$. We get $|D_{2}|\leq 13$ by Theorem
\ref{Hattingh}, Otherwise, there will be an independent set of
of size $8$ in $G[V_{>1}]$. As $G$ is a $(3,9)$-graph and
$N(v)$ is an independent set of size $\Delta(G)$, it can be deduced
that $G[D_{>2}]$ is a $(3,9-\Delta(G))$-graph. Otherwise, there is
an independent set of size $9$ in $G$, which consists of
$N(v)$ and $9-\Delta(G)$ vertices of $G[D_{>2}]$.

Now we show that $G$ is a $5$-regular graph. If $d(v)= 6,7$, or $8$, then
$G[D_{>2}]$ contains at least $6$, $5$, and $4$ vertices, respectively. But $G[D_{>2}]$ is a
$(3,3)$-graph,  a $(3,2)$-graph, or a $(3,1)$-graph $(s(3,1)=1,s(3,2)=3,s(3,3)=6)$, respectively, a contradiction. Therefore,
$d(v)\leq 5$. By $\delta(G)\geq 5$, it follows that
$d(v)=5$.

We claim that $D_{2}=13$, $D_{3}=7$ and $D_{>3}=\emptyset$. Since
$d(v)=5$ and $s(3,4)=8$, it follows that $G[D_{>2}]$ is a $(3,4)$-graph
and $|D_{>2}|\leq 7$. Easy to be calculated that $|D_{>1}|=20$
and $|D_{>2}|\leq 7$ and $|D_{2}|\leq 13$, then we deduce that
$D_{2}=13$ and $D_{>2}=7$. Now we prove $D_{>3}=\emptyset$. Suppose one vertex $w\in D_{>3}$. It illustrates that $v$, $w$ and $7$ vertices in
$D_{2}$ form an independent set. For convenience, we first illustrate some arguments.

\begin{claim}\label{claim1}
The bipartition $(X,Y)$ satisfies $|X|=7,|Y|=6$ and $D_{3}\subset N(X)$. For every nonadjacent pair
$W=\{w_{i},w_{j}\}\subset D_{3}$ there both exist a $WX$ edge and a
$WY$ edge.
\end{claim}
\begin{proof}
Clearly, it holds since $|D_{2}|=13$ and
the independent set of $G[D_{>1}]$ is at most $7$.
\end{proof}

\begin{claim}\label{claim2}
For any vertex $w\in D_{3}$ and any connected component
$(X_{i},Y_{i})$ of $(X,Y)$, either $N(w)\cap X_{i}=\emptyset$ or
$N(w)\cap Y_{i}=\emptyset$.
\end{claim}
\begin{proof}
Since $G[D_{2}\cup w]$ is a bipartite, Claim \ref{claim2} follows by Lemma \ref{Hattingh}.
\end{proof}

\begin{claim}\label{claim3}
For each nonadjacent pair $W=\{w_{i},w_{j}\}\subset D_{3}$, there is
a connected component $(X_{t},Y_{t})$ of $(X,Y)$ such that there exists
a $WX_{t}$ edge and a $WY_{t}$ edge in $G$. We call it property
$P(t)$.
\end{claim}
\begin{proof}
Clearly, it holds for Claim \ref{claim3}. Otherwise, there is
an independent set of size $8$ in $G[D_{2}\cup W]$.
\end{proof}

It satisfies that $|X|=7$ and $|Y|=6$ by Claim \ref{claim1}. Without
loss of generality, assume that $|X_{1}|=|Y_{1}|+1$ and
$|X_{i}|=|Y_{i}|$ for $i\neq 1$. Otherwise, $|X|\geq 8$ which
induces an independent set of size at least $8$. So
$D_{3}\subset N(X)$ by Claim \ref{claim1}. It furthermore follows
that $D_{3}\subset N(X_{1})$. Assume that $w\in D_{3}$ and $w\not\in
N(X_{1})$. There exists an independent set of size $8$ in
$G[D_{2}\cup w]$ by Claim \ref{claim2}.

We now discuss several cases
according to the structure of the connected components of $G[D_{2}]$.

\begin{case}
$|X|=|X_{1}|=7$ and $|Y|=|Y_{1}|=6$.
\end{case}

Because there is no $D_{3}Y_{1}$-edge by Claim \ref{claim2}, $G[Y\cup D_{3}]$ contains an independent set with of size $8$, a contradiction.

\begin{case}
$1\leq |X_{1}|\leq 2$.
\end{case}

If $|X_{1}|=1$, then $D_{3}\subset N(X_{1})$ and
$d(x)=7$, which contradicts $\Delta(G)=5$. If $|X_{1}|=2$, then at least $4$ vertices of $D_{3}$ are adjacent to one vertex
$x_{1}\in X_{1}$. Furthermore, $d(x_{1})$ is no less than $6$, which
contradicts $\Delta(G)=5$.

\begin{case}
$3\leq |X_{1}|\leq 6$.
\end{case}

Since $G[D_{3}]$ is a $(3,4)$-graph, it satisfies that $G$ contains no
$K_{3}$, an independent set of size $4$, or induced $6$-cycle by
Theorem \ref{lem2}. Suppose $G[D_{3}]$ contains a vertex with degree at least $4$, then it yields a $K_{3}$ or an independent set of
of size $4$. We have $\Delta(G[D_{3}])\leq 3$ and so there
are exactly the following eight graphs of $\mathcal{F}$ in Figure $1$.

\begin{figure}[!hbpt]
\begin{center}
\includegraphics[scale=0.7]{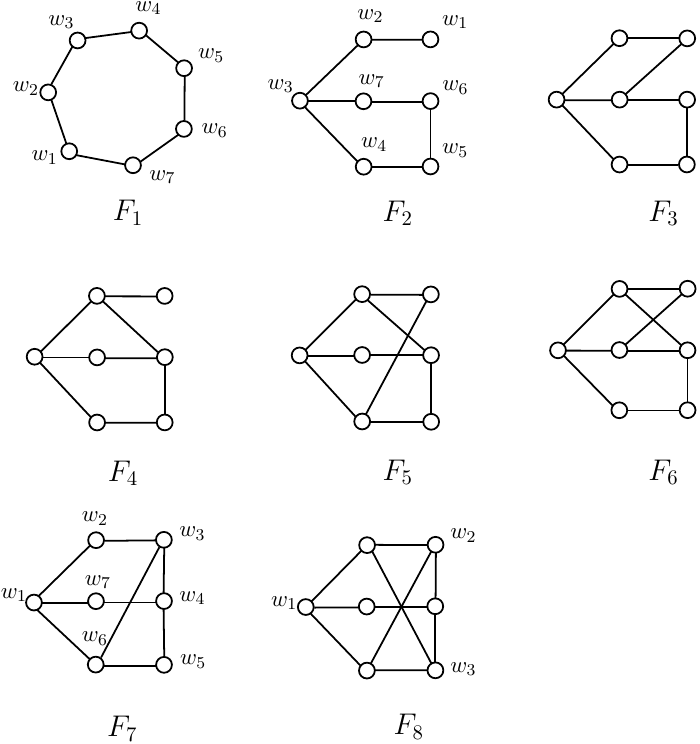}
\end{center}
\begin{center}
Figure 1: Graphs for $\mathcal{F}$.
\end{center}\label{fig4-3}
\end{figure}

Because of $D_{3}\subset N(X_{1})$, then $D_{3}\cap
N(Y_{1})=\emptyset$. Now we show that $(X,Y)$ contains at least
three connected components. Suppose, to the contrary, that the
number of connected components of $(X,Y)$ is two, say
$(X_{1},Y_{1})$ and $(X_{2},Y_{2})$. Let $\bar{D}$ denote all the
pairs of nonadjacent vertices in $D_{3}$. Thus
$(\bar{D}\cap N(X_{2}))\cup (\bar{D}\cap N(Y_{2}))=\bar{D}$ with
$\bar{D}\cap N(X_{2})\neq \emptyset$ and $\bar{D}\cap N(Y_{2})\neq
\emptyset$, which contradicts Claim \ref{claim2}.

\begin{subcase}
$G[D_{3}]=F_{1}$.
\end{subcase}

As $G[D_{3}]$ contains $14$ pairs of nonadjacent
vertices and $G[D_{2}]$ contains at most $5$ connected components,
then $3$ pairs of nonadjacent vertices in $G[D_{3}]$ have
property $P(t)$ for the same $(X_{t},Y_{t})$ by Claim \ref{claim3}.
Assume that the $3$ pairs of nonadjacent vertices in $G[D_{3}]$ are
$\{w_{1},w_{3}\},\{w_{i},w_{j}\}$ and $\{w_{r},w_{s}\}$. Then $t\neq1$ since
$D_{3}\subset N(X_{1})$. Otherwise, it contradicts
Theorem \ref{Hattingh}. Suppose $X_{t}=\{x_{1}\}$ and
$Y_{t}=\{y_{1}\}$. Without loss of generality, let
$w_{1}x_{1},w_{3}y_{1}\in E(G)$. We
have $N(w_{2})\cap X_{t}=\emptyset$ and $N(w_{2})\cap
Y_{t}=\emptyset$ since $G$ contains no triangle. Also, if $\{w_{1},w_{4}\}=\{w_{i},w_{j}\}$, then
$w_{4}y_{1}\in E(G)$, then $\{w_{3},w_{4},y_{1}\}$ is a triangle, a contradiction. Similarly,
$\{w_{3},w_{7}\}\neq\{w_{i},w_{j}\}$. Suppose
$\{w_{4},w_{7}\}=\{w_{i},w_{j}\}$. Then one vertex of $w_{5},w_{6}$
is adjacent to one vertex of $(X_{t},Y_{t})$, so it induces a
triangle. Assume that the $3$ pairs of nonadjacent
vertices in $G[D_{3}]$ are $\{w_{1},w_{4}\},\{w_{2},w_{5}\}$ and
$\{w_{3},w_{7}\}$. Without loss of generality, let
$\{y_{1}w_{2},y_{1}w_{4},y_{1}w_{7},x_{1}w_{1},x_{1}w_{3},x_{1}w_{5}\}\subset
E(G)$. Because of three pairs of nonadjacent vertices $\{w_{2},w_{4}\}$,
$\{w_{2},w_{7}\}$ and $\{w_{4},w_{7}\}$ also satisfy Claim 3
respectively, it means that one vertex of $w_{2},w_{4},w_{7}$ has
degree at least $6$, which contradicts $\Delta(G)=5$. Assume that
$\{w_{1},w_{4}\},\{w_{4},w_{6}\}$ and $\{w_{2},w_{6}\}$ satisfy
Claim 3. Let
$\{x_{1}w_{1},x_{1}w_{6},y_{1}w_{2},y_{1}w_{4}\}\subset E(G)$. If
$\{w_{3}x_{1},w_{7}y_{1}\}\cap E(G)\neq \emptyset$, then
$\{w_{1},w_{3}\},\{w_{3},w_{6}\}$ and $\{w_{1},w_{6}\}$ or
$\{w_{2},w_{4}\},\{w_{2},w_{7}\}$ and $\{w_{4},w_{7}\}$ also satisfy
Claim 3, which yields a vertex with degree at least 6 from
$w_{1},w_{3},w_{6}$ or $w_{2},w_{4},w_{7}$. If
$\{w_{3}x_{1},w_{7}y_{1}\}\cap E(G)= \emptyset$, then the pair of
vertices $\{w_{3},w_{5}\}$ satisfy Claim $3$. Furthermore, the three
pairs of vertices $\{w_{1},w_{3}\},\{w_{3},w_{6}\},\{w_{1},w_{6}\}$
also satisfy Claim \ref{claim3}, which yields a vertex of degree at
least $6$.

So it means that $|X_{t}|=|Y_{t}|\geq 2$. Suppose
$|X_{t}|=|Y_{t}|=2$ and $C=4$. Since at most two pairs of
nonadjacent vertices satisfy Claim \ref{claim3} for the same
$(X_{t},Y_{t})$ with $|X_{t}|=|Y_{t}|=1$, there are at least $10$
pairs of nonadjacent vertices in $G[D_{3}]$ connecting to the same
$(X_{t},Y_{t})$ by Claim \ref{claim3}. Without loss of generality,
if the $10$ pairs of nonadjacent vertices, $\{w_{1},w_{i}\}$ for
$i=3,4,5,6$, $\{w_{2},w_{i}\}$ for $i=4,5,6$ and $\{w_{7},w_{i}\}$
for $i=3,4,5$, satisfy Claim \ref{claim3} for the same
$(X_{t},Y_{t})$ and $w_{1}x_{1}\in E(G)$, then the pairs of
$\{w_{3},w_{5}\},\{w_{3},w_{6}\}$ are connected to one component
with of size 2 of $(X,Y)$, which will induce a triangle. Otherwise, a vertex with degree $6$ is contained in $G$. So $|(X_{i},Y_{i})|\geq 4$ for $i=1,2,\ldots,C$ and $C=3$.
Seven pairs of nonadjacent vertices, $\{w_{1},w_{i}\}$ (for
$i=3,4,5,6$), $\{w_{7}, w_{j}\}$ (for $j=3,4,5$), satisfy Claim
\ref{claim3} for the same $(X_{t},Y_{t})$. Let $X_{t}=\{x_{1},x_{2}\}$ and $Y_{t}=\{y_{1},y_{2}\}$
with
$\{x_{1}y_{1},x_{2}y_{1},x_{2}y_{2},y_{1}w_{3},y_{1}w_{5},y_{2}w_{4},y_{2}w_{6},x_{1}w_{1},x_{2}w_{7}\}\subset
E(G)$. Thus, $w_{1}x_{1}y_{1}w_{5}w_{6}w_{7}x_{1}$ is an induced
$6$-cycle. If $|X_{t}|=|Y_{t}|=3$ or $4$,
then it will contain a triangle or a vertex with degree $6$.

\begin{subcase}
$G[D_{3}]=F_{2},F_{3},F_{4},F_{5}$ or $F_{6}$.
\end{subcase}

Firstly, $G[D_{3}]$ must contain at least $12$ pairs of
nonadjacent vertices and $G[D_{2}]$ contains at most $5$ connected
components. Suppose that $C=5$, namely, Four connected
components of $G_{2}$ are both of size $2$. Since $D_{3}\subset
N(X_{1})$, $3$ pairs of nonadjacent vertices contained in $G[D_{3}]$
satisfy property $P(t)$ for $(X_{t},Y_{t})$ by Claim \ref{claim3} at
the same time. Assume that $X_{t}=\{x_{1}\},Y_{t}=\{y_{1}\}$. We can
easily check that $\{w_{3}, w_{i}\}$ for $i=1,5,6$ are three pairs
of nonadjacent vertices by Figure $1$. Without loss of generality,
let $\{w_{3}y_{1},w_{1}x_{1},w_{5}x_{1},w_{6}x_{1}\}\subset E(G)$.
Clearly, there is a triangle $w_{5}w_{6}x_{1}w_{5}$ in $G$. Thus the pair of vertices $\{w_{3},w_{5}\}$ are
connected to another connected component of $G[D_{2}]$ except
$(X_{1},Y_{1})$. Based on $D_{3}\subset N(X_{1})$,
$d(w_{3})\geq 6$, a contradiction. Then a connected
component of $G({D_{2}})$ has of size $4$ and $C\leq 4$. Let
$X_{t}=\{x_{1},x_{2}\}$ and $Y_{t}=\{y_{1},y_{2}\}$. We may let $\{w_{6}x_{2},w_{5}x_{1}\}\subset E(G)$. Clearly,
an induced 6-cycle $w_{3}w_{4}w_{5}x_{1}w_{1}w_{2}w_{3}$ occurs. It implies that the of size of a component of
$G({D_{2}})$ is 6 with $C\leq 3$. Suppose $C=3$,
then a component of $G[D_{2}]$ has of size
two. Let $X_{t}=\{x_{1},x_{2},x_{3}\}$,
$Y_{t}=\{y_{1},y_{2},y_{3}\}$. Let
$\{w_{3}y_{1},w_{1}x_{1},w_{5}x_{2},w_{6}x_{3}\}\subset E(G)$.
Therefore, the pairs of nonadjacent vertices both
$\{w_{1},w_{5}\}$ and $\{w_{1},w_{6}\}$ need to be adjacent to each
partition of the third connected component of $(X,Y)$, but then they
will induce a triangle. Suppose $C=2$. Since
$D_{3}\subset N(X_{1})$, then $D_{3}\cap
N(Y_{1})=\emptyset$. Clearly, all the pairs of nonadjacent vertices
of $D_{3}$ are connected to each partition of $(X_{t},Y_{t})$,
however, a vertex of $D_{3}$ will be adjacent to each partition of
$(X_{t},Y_{t})$, which contradicts Claim \ref{claim2}.

\begin{subcase}
$G[D_{3}]=F_{7}$ or $F_{8}$.
\end{subcase}

Suppose $G[D_{3}]=F_{7}$. Since $\{w_{1}, w_{3}\}$ satisfies
property $P(t)$ by Claim \ref{claim3}, then, without loss of
generality, let $\{w_{1}x_{1},w_{3}y_{1}\}\subset E(G)$ for
$x_{1}\in X_{t}$, $y_{1}\in Y_{t}$. Suppose
$|X_{t}|=|Y_{t}|=1$. Then for the pair of nonadjacent vertices
$\{w_{1},w_{5}\}$ it follows that $w_{5}y_{1}\notin E(G)$
(Otherwise, $K_{3}\subset G$.) and it has property $P(t')$
by Claim \ref{claim3}. Thus, $d(w_{1})\geq 6$ since $D_{3}\subset N(X_{1})$, a contradiction. Suppose
$|X_{t}|=|Y_{t}|\geq 2$. For the pair of $\{w_{1},w_{5}\}$ if
$w_{5}y_{2}\in E(G)$ for $y_{2}\in Y_{t}$, then for the pair of
$\{w_{3},w_{5}\}$ it has property $P(t'')$ by Claim \ref{claim3},
and then $d(w_{3})\geq 6$.

Suppose $G[D_{3}]=F_{8}$. Since the pairs of nonadjacent
vertices $\{w_{1},w_{2}\},\{w_{1},w_{3}\},\{w_{2},w_{3}\}$ satisfy
Claim \ref{claim3}, similarly we could deduce that there will be a
vertex of $\{w_{1},w_{2},w_{3}\}$ with degree at least $6$.
This completes the proof of Theorem \ref{thm2}.\qed

\section{A lower bound for $s_{\operatorname{CO}}(m,n)$}\label{section1}

In this section, we first prove a lower bound for $s_{\operatorname{CO}}(m,n)$ analogous to that of Krivelevich \cite{Krivelevich} for $s(m,n)$.
Let's first rephrase the definition of CO-irredundant set. Let $G=(V,E)$ be a simple graph.
Let $N[v]= \{v\}\cup N(v)$ for each $v\in V$. The open neighborhood $N(X)$ (resp., closed neighborhood $N[X]$) of a subset $X$ of $V$ is defined by $N(X) =\cup_{x\in X}N(x)$ (resp., $N[X] =\cup_{x\in X}N[x]$).
A set $X$ is called \emph{CO-irredundant} 
if, for each vertex $x\in X$, the private neighborhood of $x$ relative to $X$ satisfies that
$PN(x,X)=N[x]-N(X-\{x\})\neq \emptyset.$
Note that $y\in PN(x,X)$ if, and only if, it satisfies one of the following three cases:
\begin{itemize}
  \item[] $(i)$ $y=x$ and $x$ is an isolated vertex in $G[X]$, call $y$ a \emph{private neighbor} of $X$;
  \item[] $(ii)$ $y\in V\setminus X$ and $N(v)\cap X=\{x\}$, call $y$ a \emph{internal private neighbor} of $X$;
\item[] $(ii)$ $y\in X$ and $N(v)\cap X=\{x\}$, call $y$ an \emph{external private neighbor} of $X$.
\end{itemize}

For sets $S$ and $T$,  we denote by $K_{|S|}$ the complete graph on vertex set $S$, and denote $K_{|S|,|T|}$ the complete bipartite graph on vertex set $S$ and $T$. The \emph{completion} of $S\subset V(G)$ in $G$, denoted by $G\{S\}$, is the graph $G\cup K_{|S|}$. We denote a \emph{$k$-matching} by $kK_{2}$, and then $G-kK_{2}$ means that we delete a $k$-matching in $G$.
We first discuss that if $\langle B\rangle$ contains an $m$-element $\operatorname{CO}$-irredundant set,
then what subgraph structure $\langle R\rangle$ contains.
We now give an example of $12$-element $\operatorname{CO}$-irredundant sets in $\langle B\rangle$; see Figure 2.
\begin{example}\label{example}
Let $W_{1},W_{2},Y_{1},Y_{2},Y_{3}$ represent the vertex set $\{v_{1},v_{2},v_{3},v_{9},v_{10},v_{11},v_{12}\}$, $\{v_{5},v_{7}\}$,
$\{v_{4},v_{6},v_{8}\}$, $\{v_{4},v_{6}\}$, $\{u_{1},u_{2},u_{3},u_{4}\}$, respectively. Thus, $\langle R\rangle$ must contains the graph
$((K_{|W_{1}|}-2K_{2})\vee(K_{|Y_{1}|,|Y_{3}|}-|Y_{1}|K_{2}))\vee(K_{|W_{2}|,|Y_{2}|}\{Y_{2}\}-|W_{2}|K_{2})$.
\begin{center}
\begin{tikzpicture}
\draw [thick, blue,rounded corners] (-3,3.5) rectangle (5,1.5);
\filldraw (-2.5,2.5) circle (0.05);
\node [above] at (-2.5,2.5) {$v_{1}$};
\filldraw (-2,2.5) circle (0.05);
\node [above] at (-2,2.5) {$v_{2}$};
\filldraw (-1.5,2.5) circle (0.05);
\node [above] at (-1.5,2.5) {$v_{3}$};
\filldraw [thick] (0,2.5)circle (0.05)--(0.5,2.5)circle (0.05);
\node [above] at (0,2.5) {$v_{4}$};
\node [above] at (0.5,2.5) {$v_{5}$};
\filldraw [thick] (1,2.5)circle(0.05)--(1.5,2.5)circle(0.05);
\node [above] at (1.5,2.5) {$v_{7}$};
\node [above] at (1,2.5) {$v_{6}$};
\filldraw [thick] (1,2)circle(0.05)--(1.5,2.5)circle(0.05);
\filldraw [thick] (1,2)circle(0.05)--(0.5,2.5)circle(0.05);
\node [below,right] at (1,2) {$v_{8}$};
\filldraw [thick] (2.5,2.5)circle(0.05)--(3,2.5) circle(0.05);
\node [above] at (2.5,2.5) {$v_{9}$};
\filldraw [thick] (3.5,2.5)circle(0.05)--(4,2.5) circle(0.05);
\node [above] at (3,2.5) {$v_{10}$};
\filldraw [thick] (0,0.5)circle(0.05)--(0,2.5)circle(0.05);
\node [above] at (3.5,2.5) {$v_{11}$};
\node [above] at (4,2.5) {$v_{12}$};
\filldraw [thick] (1,0.5)circle(0.05)--(1,2)circle(0.05);
\filldraw [thick] (1.5,0.5)circle(0.05)--(1.5,2.5)circle(0.05);
\node [below] at (0,0.5){$u_{1}$};
\node [below] at (1,0.5){$u_{3}$};
\node [below] at (1.5,0.5){$u_{4}$};
\node [below] at (0.5,0.5){$u_{2}$};
\filldraw [thick](1,2.5)circle(0.05)--(0.5,0.5)circle(0.05);
\draw [thick,green,rounded corners](-0.5,1) rectangle (2,0);
\node [below] at (0.8,-1) {Figure 2: A $\operatorname{CO}$-irredundant set $\operatorname{COIR}_{12}$};
\end{tikzpicture}
\end{center}
\end{example}

So we have the following Lemma.

\begin{lemma}\label{observation1}
Let $\langle R\rangle$ and $\langle B\rangle$ be the subgraph induced by red and blue edges in a red/blue-colored complete graph, respectively. If $\langle B\rangle$ contains an $m$-element $\operatorname{CO}$-irredundant set, then $\langle R\rangle$ contains $K_{m}$ or one graph from $((K_{|W_{1}|}-tK_{2})\vee(K_{|Y_{1}|,|Y_{3}|}-|Y_{1}|K_{2}))\vee(K_{|W_{2}|,|Y_{2}|}\{Y_{2}\}-|W_{2}|K_{2})$ for $t\leq |W_{1}|\leq m-2$ as its subgraph.
\end{lemma}
\begin{proof}
Let $X=\{v_{1},v_{2},\ldots,v_{m}\}$ be an $m$-element $\operatorname{CO}$-irredundant set in $\langle B\rangle$. Clearly, if $X$ is an indpendent set in $\langle B\rangle$, then $\langle R\rangle$ contains a complete graph $K_{m}$.

We use the symbols of Figure 1. Since $W_{1}$ is the vertex-set induced by all isolated vertices and $t$ isolated edges in $\operatorname{COIR}_{m}$ of $\langle B\rangle$, there is a subgraph $K_{|W_{1}|}-tK_{2}$ for $t\leq |W_{1}|\leq m-2$ in $\langle R\rangle$. Since the subgraph induced by the edges between $Y_{1}$ and $Y_{3}$ is made up of isolated edges, then there is a subgraph $K_{|Y_{1}|,|Y_{3}|}-|Y_{1}|K_{2}$ in $\langle R\rangle$.
Let $Y_{2}$ be the internal private neighbor set. Since each pair of vertices is not connected by an edge in $Y_{2}$, then there is a $K_{|Y_{2}|}$ in $\langle R\rangle$. Then all edges between $W_{2}$ and $Y_{2}$ are isolated in $\langle B\rangle$. Moreover, there must be a subgraph $K_{|W_{2}|,|Y_{2}|}\{Y_{2}\}-|W_{2}|K_{2}$ in $\langle R\rangle$.

By the definition of $\operatorname{CO}$-irredundant set, we have that $\langle R\rangle$ contains the subgraph $((K_{|W_{1}|}-tK_{2})\vee(K_{|Y_{1}|,|Y_{3}|}-|Y_{1}|K_{2}))\vee(K_{|W_{2}|,|Y_{2}|}\{Y_{2}\}-|W_{2}|K_{2})$.
\end{proof}

For a fixed graph $H$, let $\rho(H)=\frac{|E(H)|-1}{|V(H)|-2}.$
For a collection of fixed graphs $\mathcal{H}=\{H_{1},H_{2},\ldots,H_{l}\}$, we set
$\rho(\mathcal{H})=\min\{\rho(H_{i}):1\leq i\leq l\}.$
Let $G$ be a graph. For every two disjoint subsets of $S,T\subseteq V(G)$, we denote by $e(S,T)$ the number of edges between $S$ and $T$.
\begin{theorem}[{\upshape {\bf Krivelevich} \cite{Krivelevich}}]\label{mainlemma}
Let $\mathcal{H}=\{H_{1},H_{2},\ldots,H_{k}\}$ be a family of graphs with density $\rho(\mathcal{H})>0$. Then there exists a constant $c=c(\mathcal{H})$ such that for every sufficiently large integer $N$, there exists a graph $G_{0}$ of order $N$ satisfying the following:
\begin{itemize}
\item $G_{0}$ is $\mathcal{H}$-free;

\item $G_{0}$ has no independent set of size $n=\lceil cN^{1/\rho(\mathcal{H})}\ln N\rceil$;

\item for every two disjoint subsets of $S,T\subseteq V(G)$ with $|S|=|T|=n$, we have $e(S,T)>n$.
\end{itemize}
\end{theorem}

In the following, we use Theorem \ref{mainlemma}  to prove our main theorem in this subsection.\\[0.2mm]

\noindent {\bf Proof of Theorem \ref{CO-irredund-thm}.}
For the sake of convenience, we prove there is a positive constant $c'_{m}$ such that $s_{\operatorname{CO}}(m,2n-1)>c'_{m}\left(\frac{n}{\log n}\right)^{\rho(\mathcal{H})}$ for each $m\geq 3$.

Consider an $\mathcal{H}$-free graph $G_{0}$ on $N$ vertices from Theorem \ref{mainlemma}.
We denote $G_{0}$ by $\langle R\rangle$ and denote $\overline{G_{0}}$ by $\langle B\rangle$. Then we have the following.

\begin{claim}\label{claim4}
$\langle B\rangle$ does not contain an $m$-element $\operatorname{CO}$-irredundant set.
\end{claim}
\begin{proof}
Assume, to the contrary, that $\langle B\rangle$ contains an $\operatorname{COIR}_m$. By Lemma \ref{observation1}, $\langle R\rangle$ contains $K_{m}$ or one graph from
$\{((K_{|W_{1}|}-tK_{2})\vee(K_{|Y_{1}|,|Y_{3}|}-|Y_{1}|K_{2}))\vee(K_{|W_{2}|,|Y_{2}|}\{Y_{2}\}-|W_{2}|K_{2})|\ t\leq |W_{1}|\leq m-2\}$. We apply Theorem \ref{mainlemma} to $\mathcal{H}=\{((K_{|W_{1}|}-tK_{2})\vee(K_{|Y_{1}|,|Y_{3}|}-|Y_{1}|K_{2}))\vee(K_{|W_{2}|,|Y_{2}|}\{Y_{2}\}-|W_{2}|K_{2})|\ t\leq |W_{1}|\leq m-2\}$, a contradiction.
\end{proof}

\begin{claim}\label{claim5}
$\langle R\rangle$ does not contain an $2n-1$-element $\operatorname{CO}$-irredundant set.
\end{claim}
\begin{proof}
Assume, to the contrary, that $\langle R\rangle$ contains a $\operatorname{CO}$-irredundant set, say $D$, of size $2n-1$.

By Lemma \ref{mainlemma}, there exists an independent set, say $D'$, of size at most $n-2$. Fix the independent set $D'$ of size $n-2$. For the remaining $n+1$ vertices, there is at most one leaf in the subgraph induced by the remaining $n+1$ vertices. Otherwise, $\langle R\rangle$ contains an independent set of size $n$, a contradiction.

It means that there is at most one vertex containing an internal private neighbor. Let $v$ be a vertex adjacent to the leaf. each of the remaining $n$ vertices (except $v$), say $D'$, must be adjacent to an external private neighbor. We denote by $W$ the set of external private neighbors. Therefore, $e(D',W)=n$, which contradicts Theorem \ref{mainlemma}.
\end{proof}

From Claims \ref{claim4} and \ref{claim5}, $\langle B\rangle$ does not contain a $\operatorname{CO}$-irredundant set of size $m$ or $\langle R\rangle$ does  a CO-irredundant set of size $2n-1$.\qed
\vskip 0.5cm
\noindent{\bf Acknowledgment.} We are also extremely grateful to Prof. Xian'an Jin for carefully reading an earlier draft of this paper.


\end{document}